\documentclass[12pt,twoside]{amsart}
\usepackage[latin1]{inputenc}
\usepackage{amsmath, amsthm, amscd, amsfonts, amssymb, graphicx}
\usepackage[bookmarksnumbered, plainpages]{hyperref}

\newtheorem{thm}{Theorem}[section]

\newtheorem{defn}[thm]{Definition}

\numberwithin{equation}{section}

\begin{document}

\title{\bf Geometric structures of Morris-Thorne wormhole metric with the semi-symmetric non-metric connections}
\author{Siyao Liu  \hskip 0.4 true cm  Yong Wang$^{*}$}

\thanks{{\scriptsize
\hskip -0.4 true cm \textit{2010 Mathematics Subject Classification:}
53B20; 53B25; 53B30; 53C15; 53C25; 53C35.
\newline \textit{$^{*}$Corresponding author}}}

\maketitle

\begin{abstract}
 \indent Spacetime is a $4$-dimensional connected Lorentzian manifold. In this paper, we extend the Levi-Civita connection in the definition of spacetime to the semi-symmetric non-metric connection and conclude geometric structures admitted by the metric $(1.1)$ with the semi-symmetric non-metric connections. Obviously, Morris-Thorne spacetime is Ricci generalized pseudosymmetry, Ricci generalized projectively pseudosymmetry, and has conformal $2$-forms that are recurrent, etc. It also is an Einstein manifold of level $2$ and $3$-quasi-Einstein manifold.
\end{abstract}

\vskip 0.2 true cm

%------------------------------------------------------------------------------------%

\pagestyle{myheadings}
\markboth{\rightline {\scriptsize  Liu}}
         {\leftline{\scriptsize }}

\bigskip
\bigskip

%------------------------------------------------------------------------------------%
%------------------------------------------------------------------------------------%

\section{ Introduction}
Einstein and Rosen defined the Einstein-Rosen bridge in $1935,$ which can connect two different points in spacetime, thus theoretically making shortcuts that can reduce travel time and distance \cite{ER}.
The notion of a wormhole was first introduced in $1962$ by Wheeler and reinterpreted the Einstein-Rosen bridge as a connection between two distant places in spacetime with no more mutual interaction \cite{W}.
With the development of the research, Wheeler and Fuller found that this Schwarzschild wormhole cannot be traversed even by a single particle \cite{FW}.
In $1973,$ Ellis first mentioned the following metric called a drain hole, or a topological hole \cite{E}:
\begin{align}
ds^{2}=-c^{2}dt^{2}+dl^{2}+(b^{2}+l^{2})dv^{2}+(b^{2}+l^{2})sin^{2}vd\phi^{2},
\end{align}
where $t$ is the global time, $l$ is the proper radial coordinate, $b$ is the shape constant and $c$ is the speed of light.
Morris and Thorne recognized metric $(1.1)$ as a wormhole that was traversed in principle by human beings \cite{MT}.
Based on the importance of the Morris-Thorne wormhole metric in astrophysics, the researchers have given many conclusions related to it \cite{LLO,T1,J}.

It is only natural that the researchers investigate the geometric properties of the Morris-Thorne wormhole metric, that is describe the geometric structures of corresponding spacetimes, such as semisymmetry \cite{S1,S2,S3}, Deszcz pseudosymmetry \cite{AD}, Chaki pseudosymmetry \cite{C} and recurrent manifold \cite{R1,R2,R3}.
Sabina in \cite{SDM} gave the geometric properties admitted by Morris-Thorne spacetime, such as Ricci generalized pseudosymmetry, Ricci generalized projectively pseudosymmetry, pseudosymmetry due to Weyl conformal curvature, semisymmetry due to conharmonic curvature, etc.

Agashe and Chafle introduced the notion of a semi-symmetric non-metric connection and studied some of its properties and submanifolds of a Riemannian manifold with a semi-symmetric non-metric connection \cite{AC1,AC2}.
Wang also researched non-integrable distributions with semi-symmetric non-metric connections \cite{Wang1}. 
In \cite{Wang2}, the author considered the generalized Kasner spacetimes with a semi-symmetric non-metric connection.

Inspired by the work of Sabina, Agashe, Chafle and Wang, we extend the Levi-Civita connection in the definition of spacetime to the semi-symmetric non-metric connection and consider the geometric properties in terms of curvatures admitted by this spacetime.
The method of this paper is to deduce the different components of various tensors of the metric $(1.1)$ and conclude the curvature restricted geometric structures admitted by the metric $(1.1).$
A trivial verification shows that the metric $(1.1)$ realizes several important geometric structures, for instance, it is Ricci generalized pseudosymmetry, Ricci generalized projectively pseudosymmetry and it has conformal $2$-forms that are recurrent, etc.
It also is an Einstein manifold of level $2$ and $3$-quasi-Einstein manifold.

A brief description of the organization of this paper is as follows.
We explain in Sec. 2 the basic notions of several geometric structures which need to be used in this paper.
In Section 3, we provide a detailed exposition of the computational process of various tensors and thus prove the main results.
Section 4 contains a summary of the main conclusions of this article.

\vskip 1 true cm

\section{Preliminaries}
Let $M$ be a connected smooth semi-Riemannian manifold of dimension $n(\geq 3)$ equipped with the semi-Riemannian metric $g.$
From now on, $\nabla$ denotes the Levi-Civita connection and $X, Y, Y_{1}, Y_{2}..., U, U_{1}, U_{2}..., V, V_{1}, V_{2}..., \widehat{P}\in \mathfrak{X}(M).$
We consider the semi-symmetric non-metric connection on $M$ \cite{Wang1,Wang2}
\begin{align}
\widehat{\nabla}_{X}Y=\nabla_{X}Y+g(\widehat{P},Y)X.
\end{align}
$A$ and $E$ are two symmetric $(0,2)$-tensors, we define their Kulkarni-Nomizu product $A\wedge E$ as \cite{G1,G2,SRK}
\begin{align}
(A\wedge E)(Y_{1}, Y_{2}, U_{1}, U_{2})&=A(Y_{1}, U_{2})E(Y_{2},U_{1})-A(Y_{1},U_{1})E(Y_{2}, U_{2})\\
&+A(Y_{2},U_{1})E(Y_{1}, U_{2})-A(Y_{2}, U_{2})E(Y_{1}, U_{1}).\nonumber
\end{align}
For a symmetric $(0,2)$-tensor $A,$ we can define the endomorphism $Y_{1}\wedge_{A} Y_{2}$ as \cite{DDHKS,SH1,SH2}
\begin{align}
(Y_{1}\wedge_{A} Y_{2})U=A(Y_{2}, U)Y_{1}-A(Y_{1},U)Y_{2}.
\end{align}
Write
\begin{align}
\mathcal{E}_{R}&(Y_{1}, Y_{2})=\widehat{\nabla}_{Y_{1}}\widehat{\nabla}_{Y_{2}}-\widehat{\nabla}_{Y_{2}}\widehat{\nabla}_{Y_{1}}-\widehat{\nabla}_{[Y_{1}, Y_{2}]}.
\end{align}
We introduce four endomorphisms on $M$ \cite{K,SC}
\begin{align}
\mathcal{E}_{C}&=\mathcal{E}_{R}-\frac{1}{n-2}(\mathcal{J}\wedge_{g}+\wedge_{g}\mathcal{J}-\frac{\kappa}{n-1}\wedge_{g}),\\
\mathcal{E}_{K}&=\mathcal{E}_{R}-\frac{1}{n-2}(\mathcal{J}\wedge_{g}+\wedge_{g}\mathcal{J}),\\
\mathcal{E}_{W}&=\mathcal{E}_{R}-\frac{\kappa}{n(n-1)}\wedge_{g},\\
\mathcal{E}_{P}&=\mathcal{E}_{R}-\frac{1}{n-1}\wedge_{Ric},
\end{align}
where $\mathcal{J}$ is the Ricci operator defined as $Ric(Y_{1}, Y_{2})=g(Y_{1}, \mathcal{J}Y_{2}),$ 
the Ricci tensor $Ric$ is defined by $Ric(Y_{1}, Y_{2})=tr\{X\rightarrow \mathcal{E}_{R}(X, Y_{1})Y_{2}\}$
and $\kappa$ denotes the scalar curvature with the semi-symmetric non-metric connections.
For an endomorphism $\mathcal{E}(U_{1}, U_{2})$ one can define the following $(0,4)$ type tensor field as
\begin{align}
E(U_{1}, U_{2}, U_{3}, U_{4})=g(\mathcal{E}(U_{1}, U_{2})U_{3}, U_{4}).
\end{align}
If we replace $\mathcal{E}$ with the endomorphisms $\mathcal{E}_{R}$ $\mathcal{E}_{C}$ $\mathcal{E}_{K}$ $\mathcal{E}_{W}$ and $\mathcal{E}_{P}$ above, then we can  obtain the $(0,4)$ type Riemann curvature $R,$ the Weyl conformal curvature $C,$ the conharmonic curvature $K,$ the concircular curvature $W$ and the projective curvature $P$ about $\widehat{\nabla}.$
These tensors locally are given by
\begin{align}
R_{hkij}&=g_{h\alpha}\left(\partial_{i}\widehat{\Gamma}^{\alpha}_{kj}-\partial_{j}\widehat{\Gamma}^{\alpha}_{ki}+\widehat{\Gamma}^{\beta}_{kj}\widehat{\Gamma}^{\alpha}_{\beta i}-\widehat{\Gamma}^{\beta}_{ki}\widehat{\Gamma}^{\alpha}_{\beta j}\right),\\
C_{hkij}&=\left(R-\frac{1}{n-2}(g\wedge Ric)+\frac{\kappa}{2(n-1)(n-2)}(g\wedge g)\right)_{hkij},\\
K_{hkij}&=\left(R-\frac{1}{n-2}(g\wedge Ric)\right)_{hkij},\\
W_{hkij}&=\left(R-\frac{\kappa}{2n(n-1)}(g\wedge g)\right)_{hkij},\\
P_{hkij}&=R_{hkij}-\frac{1}{n-1}(g_{hi}Ric_{kj}-g_{ki}Ric_{hj}),
\end{align}
where $g_{h\alpha}$ is the metric matrix, $\partial_{i}$ is a natural local frame, $\widehat{\Gamma}^{\alpha}_{kj}$ is the Christofel coefficient of the semi-symmetric non-metric connection.
Here $F$ is a $(0,k),$ $k\geq 1,$ type tensor field, we have a $(0,k+2)$ type tensor field $E\cdot F$ as follows \cite{DG,SH3}:
\begin{align}
(E\cdot F)(Y_{1}, Y_{2},..., Y_{k}, U_{1}, U_{2})&=(\mathcal{E}(U_{1}, U_{2})\cdot F)(Y_{1}, Y_{2},..., Y_{k})\\
&=-F(\mathcal{E}(U_{1}, U_{2})Y_{1}, Y_{2},..., Y_{k})-...-F(Y_{1}, Y_{2},..., \mathcal{E}(U_{1}, U_{2})Y_{k}),\nonumber
\end{align}
and a $(0,k+2)$ type tensor $Q(E,F),$ named Tachibana tensor, as follows \cite{T2,DGPSS}:
\begin{align}
&Q(E, F)(Y_{1}, Y_{2},..., Y_{k}, U_{1}, U_{2})=((U_{1}\wedge_{E} U_{2})\cdot F)(Y_{1}, Y_{2},..., Y_{k})\\
&=E(U_{1}, Y_{1})F(U_{2}, Y_{2},..., Y_{k})+...+E(U_{1}, Y_{k})F(Y_{1}, Y_{2},..., U_{2})\nonumber\\
&-E(U_{2}, Y_{1})F(U_{1}, Y_{2},..., Y_{k})-...-E(U_{2}, Y_{k})F(Y_{1}, Y_{2},..., U_{1}).\nonumber
\end{align}
\begin{defn} \cite{R4,R5,SH4,SH5}
If the tensor $E\cdot F$ is linearly dependent with $Q(Z,F)$ then it is called as $F$-pseudosymmetric type manifold due to $E,$  if $E\cdot F=0$ then it is called $F$-semisymmetric type manifold due to $E.$
\end{defn}
In particular, if $E\cdot F= f_{F}Q(Ric,F)$ ($f_{F}$ being some smooth scalar function on $M$), then it is called a Ricci generalized $F$-pseudosymmetric type manifold due to $E.$
A $F$-pseudosymmetric manifold due to $E$ is said to be a Deszcz pseudosymmetric manifold if $E=R,$ $F=R$ and $Z=g,$ and a $F$-semisymmetric type manifold due to $E$ is said to be a semisymmetric manifold if $E=R,$ $F=R.$
Again for $E=C,$ $F=C$ and $Z=g$ the manifold is said to have a pseudosymmetric Weyl conformal curvature tensor.
With the definition in \cite{S,SYH}, we have
\begin{defn}
For a scalar $\alpha$ and $0\leq k\leq n-1$, if $rank(Ric-\alpha g)=k,$ then $M$ is called as a $k$-quasi-Einstein manifold and for $k=1$ (resp., $k=0$) it is a quasi-Einstein (resp., Einstein) manifold. A Ricci simple manifold is a quasi-Einstein manifold for $\alpha = 0$.
\end{defn}
Our definition is similar to the one given in \cite{B,D,DGJZ}. 
\begin{defn}
A generalized Roter type manifold $M$ is a manifold whose Riemann curvature can explicitly be expressed as
\begin{align}
R=(\mu_{11}Ric^{2}+\mu_{12}Ric+\mu_{13}g)\wedge Ric^{2}+(\mu_{22}Ric+\mu_{23}g)\wedge Ric+\mu_{33}(g\wedge g),
\end{align}
where $\mu_{ij}$ are some scalars, the Ricci tensor of level $k$ is defined by $Ric^{k}(Y_{1}, Y_{2})=g(Y_{1}, \mathcal{J}^{k-1}Y_{2}).$ If such linear dependency reduces to the linear dependency of $R,$ $g\wedge g,$ $g\wedge Ric,$ $Ric\wedge Ric$ then it is called a manifold of Roter.
\end{defn}
Since the definition in \cite{SH6}, Einstein manifold of level $2$ (resp., $3$ and $4$) about the connection $\widehat{\nabla}$ is well defined.
\begin{defn} \cite{SH6}
A semi-Riemannian manifold $M$ is said to be Einstein manifold of level $2$ (resp., $3$ and $4$) about the connection $\widehat{\nabla}$ if
\begin{align}
&Ric^{2}+\lambda_{1}Ric+\lambda_{2}g=0,\\
&Ric^{3}+\lambda_{3}Ric^{2}+\lambda_{4}Ric+\lambda_{5}g=0,\\
&Ric^{4}+\lambda_{6}Ric^{3}+\lambda_{7}Ric^{2}+\lambda_{8}Ric+\lambda_{9}g=0,
\end{align}
where $\lambda_{i}(1\leq\lambda_{i}\leq9)$ are some scalars.
\end{defn}
Based on the concept of \cite{G,SJ}, we have
\begin{defn}
The Ricci tensor of the semi-Riemannian manifold $M$ is said to be cyclic parallel if the relation
\begin{align}
\underset{Y_{1}, Y_{2}, Y_{3}}{\mathcal{S}} (\widehat{\nabla}_{Y_{1}}Ric)(Y_{2}, Y_{3})=0
\end{align}
holds and it is said to be Codazzi if $(\widehat{\nabla}_{Y_{1}}Ric)(Y_{2}, Y_{3})=(\widehat{\nabla}_{Y_{2}}Ric)(Y_{1}, Y_{3})$ holds where $\mathcal{S}$ is a cyclic sum over $Y_{1},$ $Y_{2}$ and
$Y_{3}.$
\end{defn}
\begin{defn} \cite{MM1,MM2,MM3}
A symmetric $(0,2)$-tensor $Z$ on a semi-Riemannian M is said to be $T$-compatible if
\begin{align}
\underset{Y_{1}, Y_{2}, Y_{3}}{\mathcal{S}} T(\mathcal{Z}Y_{1}, U, Y_{2}, Y_{3})=0
\end{align}
holds, where $\mathcal{Z}$ is the endomorphism corresponding to $Z$ defined as $g(\mathcal{Z}Y_{1}, Y_{2})= Z(Y_{1}, Y_{2}).$
\end{defn}
Replacing $T$ by the curvatures $R,$ $C,$ $W,$ $P$ and $K$ we can define the corresponding curvature compatibilities.
We develop the definition of \cite{MS1,MS2}, then  we could find
\begin{defn}
The curvature $2$-forms $\Omega^{M}_{(E)l}$ are recurrent if and only if
\begin{align}
\underset{Y_{1}, Y_{2}, Y_{3}}{\mathcal{S}} (\widehat{\nabla}_{Y_{1}}E)(Y_{2}, Y_{3}, U, Y)=\underset{Y_{1}, Y_{2}, Y_{3}}{\mathcal{S}} \sigma(Y_{1})E(Y_{2}, Y_{3}, U, Y).
\end{align}
\end{defn}

\vskip 1 true cm

\section{Various curvature tensors of Morris-Thorne wormhole metric}
We will write the spherical coordinates $(t, l, v, \phi)$ as $(X_{1}, X_{2}, X_{3}, X_{4})$ and let $\widehat{P}=a\frac{\partial}{\partial X_{2}}=a\partial_{2},$ where $a$ is a constant.
The components of the Morris-Thorne metric $(1.1)$ are given by
\begin{align}
g_{11}=-c^{2};~g_{22}=1;~g_{33}=b^{2}+X_{2}^{2};~g_{44}=(b^{2}+X_{2}^{2})\sin^{2}X_{3};~g_{ij}=0,otherwise.
\end{align}

The non-vanishing components of the Christofel coefficient $\widehat{\Gamma}^{\alpha}_{kj}$ are given below
\begin{align}
&\widehat{\Gamma}^{1}_{12}=\widehat{\Gamma}^{2}_{22}=a,~\widehat{\Gamma}^{2}_{33}=-X_{2},~\widehat{\Gamma}^{4}_{34}=\widehat{\Gamma}^{4}_{43}=\cot X_{3},\\
&\widehat{\Gamma}^{3}_{23}=\widehat{\Gamma}^{4}_{24}=\frac{X_{2}}{b^{2}+X_{2}^{2}},~\widehat{\Gamma}^{3}_{32}=\widehat{\Gamma}^{4}_{42}=\frac{X_{2}}{b^{2}+X_{2}^{2}}+a,\nonumber\\
&\widehat{\Gamma}^{2}_{44}=-X_{2}\sin^{2}X_{3},~\widehat{\Gamma}^{3}_{44}=-\sin X_{3}\cos X_{3}.\nonumber
\end{align}
Then the non-zero components of the Riemann curvature tensor $R,$ Ricci tensor $Ric$ and scalar curvature $\kappa$ of (1.1) are given by
\begin{align}
&R_{2323}=-\frac{b^{2}}{b^{2}+X_{2}^{2}},~R_{2424}=-\frac{b^{2}\sin^{2}X_{3}}{b^{2}+X_{2}^{2}},~R_{3434}=b^{2}\sin^{2}X_{3};\\
&Ric_{22}=\frac{2b^{2}}{(b^{2}+X_{2}^{2})^{2}};~\kappa=\frac{2b^{2}}{(b^{2}+X_{2}^{2})^{2}}.\nonumber
\end{align}
Thus, the non-zero components of $\widehat{\nabla} R$ and $\widehat{\nabla} Ric$ are calculated as below
\begin{align}
&\widehat{\nabla}_{2} R_{2323}=\frac{2b^{2}(ab^{2}+2X_{2}+aX_{2}^{2})}{(b^{2}+X_{2}^{2})^{2}},~\widehat{\nabla}_{2} R_{2424}=\frac{2b^{2}\sin^{2}X_{3}(ab^{2}+2X_{2}+aX_{2}^{2})}{(b^{2}+X_{2}^{2})^{2}},\\
&\widehat{\nabla}_{2} R_{3434}=-\frac{4b^{2}X_{2}\sin^{2}X_{3}}{b^{2}+X_{2}^{2}},~\widehat{\nabla}_{3} R_{2434}=-\widehat{\nabla}_{4} R_{2334}==-\frac{b^{2}\sin^{2}X_{3}(ab^{2}+2X_{2}+aX_{2}^{2})}{b^{2}+X_{2}^{2}};\nonumber\\
&\widehat{\nabla}_{2} Ric_{22}=-\frac{4b^{2}(ab^{2}+2X_{2}+aX_{2}^{2})}{(b^{2}+X_{2}^{2})^{3}},~\widehat{\nabla}_{3} Ric_{23}=\frac{2b^{2}X_{2}}{(b^{2}+X_{2}^{2})^{2}},~\widehat{\nabla}_{4} Ric_{24}=\frac{2b^{2}X_{2}\sin^{2}X_{3}}{(b^{2}+X_{2}^{2})^{2}}.\nonumber
\end{align}
According to some computations about $g\wedge Ric$ and $g\wedge g$, we can get the non zero components of the tensor $C,$ $P,$ $K$ as follows
\begin{align}
&C_{1212}=-\frac{2b^{2}c^{2}}{3(b^{2}+X_{2}^{2})^{2}},~C_{1313}=\frac{b^{2}c^{2}}{3(b^{2}+X_{2}^{2})},~C_{1414}=\frac{b^{2}c^{2}\sin^{2}X_{3}}{3(b^{2}+X_{2}^{2})},\\
&C_{2323}=-\frac{b^{2}}{3(b^{2}+X_{2}^{2})},~C_{2424}=-\frac{b^{2}\sin^{2}X_{3}}{3(b^{2}+X_{2}^{2})},~C_{3434}=\frac{2}{3}b^{2}\sin^{2}X_{3};\nonumber\\
&P_{1221}=\frac{2b^{2}c^{2}}{3(b^{2}+X_{2}^{2})^{2}},~P_{2323}=-\frac{b^{2}}{3(b^{2}+X_{2}^{2})},~P_{2332}=\frac{b^{2}}{b^{2}+X_{2}^{2}},~\nonumber\\
&P_{2424}=-\frac{b^{2}\sin^{2}X_{3}}{3(b^{2}+X_{2}^{2})},~P_{2442}=\frac{b^{2}\sin^{2}X_{3}}{b^{2}+X_{2}^{2}},~P_{3434}=-P_{3443}=b^{2}\sin^{2}X_{3};\nonumber\\
&K_{1212}=-\frac{b^{2}c^{2}}{(b^{2}+X_{2}^{2})^{2}},~K_{3434}=b^{2}\sin^{2}X_{3}.\nonumber
\end{align}
It is easy to check that
\begin{align}
&\widehat{\nabla}_{2} C_{1212}=\frac{4b^{2}c^{2}(ab^{2}+2X_{2}+aX_{2}^{2})}{3(b^{2}+X_{2}^{2})^{3}},~\widehat{\nabla}_{3} C_{1213}=-\frac{b^{2}c^{2}(ab^{2}+3X_{2}+aX_{2}^{2})}{3(b^{2}+X_{2}^{2})^{2}},\\
&\widehat{\nabla}_{4} C_{1214}=-\frac{b^{2}c^{2}\sin^{2}X_{3}(ab^{2}+3X_{2}+aX_{2}^{2})}{3(b^{2}+X_{2}^{2})^{2}},~\widehat{\nabla}_{2} C_{1313}=-\frac{4b^{2}c^{2}X_{2}}{3(b^{2}+X_{2}^{2})^{2}},\nonumber\\
&\widehat{\nabla}_{1} C_{1323}=-\frac{ab^{2}c^{2}}{3(b^{2}+X_{2}^{2})},~\widehat{\nabla}_{2} C_{1414}=-\frac{4b^{2}c^{2}X_{2}\sin^{2}X_{3}}{3(b^{2}+X_{2}^{2})^{2}},~\widehat{\nabla}_{1} C_{1424}=-\frac{ab^{2}c^{2}\sin^{2}X_{3}}{3(b^{2}+X_{2}^{2})};\nonumber\\
&\widehat{\nabla}_{2} P_{1221}=-\frac{4b^{2}c^{2}(ab^{2}+2X_{2}+aX_{2}^{2})}{3(b^{2}+X_{2}^{2})^{3}},~\widehat{\nabla}_{1} P_{1222}=-\frac{2ab^{2}c^{2}}{3(b^{2}+X_{2}^{2})^{2}},\nonumber\\
&\widehat{\nabla}_{3} P_{1231}=\widehat{\nabla}_{3} P_{1321}=\frac{2b^{2}c^{2}X_{2}}{3(b^{2}+X_{2}^{2})^{2}},~\widehat{\nabla}_{4} P_{1241}=\widehat{\nabla}_{4} P_{1421}=\frac{2b^{2}c^{2}X_{2}\sin^{2}X_{3}}{3(b^{2}+X_{2}^{2})^{2}};\nonumber\\
&\widehat{\nabla}_{2} K_{1212}=\frac{2b^{2}c^{2}(ab^{2}+2X_{2}+aX_{2}^{2})}{(b^{2}+X_{2}^{2})^{3}},~\widehat{\nabla}_{3} K_{1213}=-\frac{b^{2}c^{2}X_{2}}{(b^{2}+X_{2}^{2})^{2}},~\widehat{\nabla}_{4} K_{1214}=-\frac{b^{2}c^{2}X_{2}\sin^{2}X_{3}}{(b^{2}+X_{2}^{2})^{2}};\nonumber
\end{align}
Again the non-zero components of the tensor $R\cdot R,$ $P\cdot R$ are written below
\begin{align}
&(R\cdot R)_{233424}=-\frac{2b^{4}\sin^{2}X_{3}}{(b^{2}+X_{2}^{2})^{2}},~(R\cdot R)_{243423}=\frac{2b^{4}\sin^{2}X_{3}}{(b^{2}+X_{2}^{2})^{2}};\\
&(P\cdot R)_{233424}=-\frac{4b^{4}\sin^{2}X_{3}}{3(b^{2}+X_{2}^{2})^{2}},~(P\cdot R)_{243423}=\frac{4b^{4}\sin^{2}X_{3}}{3(b^{2}+X_{2}^{2})^{2}};\nonumber
\end{align}
An easy computation shows that
\begin{align}
&(C\cdot C)_{121323}=-(C\cdot C)_{122313}=\frac{b^{4}c^{2}}{3(b^{2}+X_{2}^{2})^{3}},~(C\cdot C)_{121424}=-(C\cdot C)_{122414}=\frac{b^{4}c^{2}\sin^{2}X_{3}}{3(b^{2}+X_{2}^{2})^{3}},\\
&(C\cdot C)_{133414}=-(C\cdot C)_{143413}=\frac{b^{4}c^{2}\sin^{2}X_{3}}{3(b^{2}+X_{2}^{2})^{2}},~(C\cdot C)_{233424}=-(C\cdot C)_{243423}=-\frac{b^{4}\sin^{2}X_{3}}{3(b^{2}+X_{2}^{2})^{2}};\nonumber\\
&(C\cdot K)_{121323}=-(C\cdot K)_{122313}=(C\cdot K)_{131223}=-(C\cdot K)_{231213}=\frac{b^{4}c^{2}}{3(b^{2}+X_{2}^{2})^{3}},\nonumber\\
&(C\cdot K)_{121424}=-(C\cdot K)_{122414}=(C\cdot K)_{141224}=-(C\cdot K)_{241214}=\frac{b^{4}c^{2}\sin^{2}X_{3}}{3(b^{2}+X_{2}^{2})^{3}},\nonumber\\
&(C\cdot K)_{133414}=-(C\cdot K)_{143413}=\frac{b^{4}c^{2}\sin^{2}X_{3}}{3(b^{2}+X_{2}^{2})^{2}},~(C\cdot K)_{233424}=-(C\cdot K)_{243423}=-\frac{b^{4}\sin^{2}X_{3}}{3(b^{2}+X_{2}^{2})^{2}};\nonumber
\end{align}
By $(2.16),$ the non-vanishing components of the Tachibana tensor are given by
\begin{align}
&Q(Ric, R)_{233424}=-\frac{2b^{4}\sin^{2}X_{3}}{(b^{2}+X_{2}^{2})^{2}},~Q(Ric, R)_{243423}=\frac{2b^{4}\sin^{2}X_{3}}{(b^{2}+X_{2}^{2})^{2}};\\
&Q(g, C)_{121323}=-Q(g, C)_{122313}=\frac{b^{2}c^{2}}{b^{2}+X_{2}^{2}},~Q(g, C)_{121424}=-Q(g, C)_{122414}=\frac{b^{2}c^{2}\sin^{2}X_{3}}{b^{2}+X_{2}^{2}},\nonumber\\
&Q(g, C)_{133414}=-Q(g, C)_{143413}=b^{2}c^{2}\sin^{2}X_{3},~Q(g, C)_{233424}=-Q(g, C)_{243423}=-b^{2}\sin^{2}X_{3};\nonumber\\
&Q(g, K)_{121323}=-Q(g, K)_{122313}=Q(g, K)_{131223}=-Q(g, K)_{231213}=\frac{b^{2}c^{2}}{b^{2}+X_{2}^{2}},\nonumber\\
&Q(g, K)_{121424}=-Q(g, K)_{122414}=Q(g, K)_{141224}=-Q(g, K)_{241214}=\frac{b^{2}c^{2}\sin^{2}X_{3}}{b^{2}+X_{2}^{2}},\nonumber\\
&Q(g, K)_{133414}=-Q(g, K)_{143413}=b^{2}c^{2}\sin^{2}X_{3},~Q(g, K)_{233424}=-Q(g, K)_{243423}=-b^{2}\sin^{2}X_{3};\nonumber
\end{align}

\vskip 1 true cm

\section{Conclusions}
Now let us give the the main results in this paper.
\begin{thm}
The Morris-Thorne wormhole metric $(1.1)$ fulfills the following geometric structures about the connection $\widehat{\nabla}$:\\
(I)(1)Ricci generalized pseudosymmetry  i.e. $R\cdot R=Q(Ric, R);$\\
(2)Ricci generalized projectively pseudosymmetric because $P\cdot R=\frac{2}{3}Q(Ric, R);$\\
(3)Pseudosymmetric Weyl conformal curvature since $C\cdot C=\frac{b^{2}}{3(b^{2}+X_{2}^{2})^{2}}Q(g, C);$\\
(4)Conharmonic curvature pseudosymmetric type manifold due to the Weyl conformal curvature satisfying $C\cdot K=\frac{b^{2}}{3(b^{2}+X_{2}^{2})^{2}}Q(g, K);$\\
(5)$K\cdot C=0$ hence it has the Weyl conformal curvature semisymmetric type manifold due to conharmonic curvature;\\
(6)Semisymmetric conharmonic curvature as it fulfills $K\cdot K=0.$\\
(II)$Ric=\alpha(\eta\bigotimes\eta)$ for $\alpha=\frac{2b^{2}}{(b^{2}+X_{2}^{2})^2}$ and $\eta=\{0, 1, 0, 0\},$ thus it is a Ricci simple manifold and a $3$-quasi-Einstein manifold as $rank(Ric-\alpha g)=3.$\\
(III)The manifold neither Roter type nor the generalized Roter type, but $g\wedge Ric^{2}=\frac{2b^{2}}{b^{2}+X_{2}^{2}}g\wedge Ric.$\\
(IV)The condition $Ric^{2}=\frac{2b^{2}}{(b^{2}+X_{2}^{2})^{2}}Ric$ holds, hence it is a Einstein manifold of level $2.$\\
(V)Ricci tensor is neither of Codazzi type nor cyclic parallel.\\
(VI)Ricci tensor is Riemann compatible, Weyl conformal compatible and conharmonic compatible.\\
(VII)Conformal curvature $2$-forms are recurrent for the $1$-form $\{0, -\frac{X_{2}}{b^{2}+X_{2}^{2}}, 0, 0\}.$\\

\end{thm}

\vskip 1 true cm

\section{Acknowledgements}

The author was supported in part by  NSFC No.11771070. The author thanks the referee for his (or her) careful reading and helpful comments.

\vskip 1 true cm

%-----------------------------------------------------------------------------
%-----------------------------------------------------------------------------

\bigskip
\bigskip

\noindent {\footnotesize {\it S. Liu} \\
{School of Mathematics and Statistics, Northeast Normal University, Changchun 130024, China}\\
{Email: liusy719@nenu.edu.cn}

\noindent {\footnotesize {\it Y. Wang} \\
{School of Mathematics and Statistics, Northeast Normal University, Changchun 130024, China}\\
{Email: wangy581@nenu.edu.cn}


\begin{thebibliography}{00}
\bibitem{ER} A. Einstein, N. Rosen, The particle problem in the general theory of relativity, Phys. Rev. 48 (1935) 73.
\bibitem{W} A.J. Wheeler, Geometrodynamics, Academic Press, New York, 1962.
\bibitem{FW} W.R. Fuller, A.J. Wheeler, Causality and multiply-connected space-time, Phys. Rev. 128 (1962) 919.
\bibitem{E} H.G. Ellis, Ether flow through a drainhole: A particle model in general relativity, J. Math. Phys. 14 (1973) 104.
\bibitem{MT} S.M. Morris, S.K. Thorne, Wormholes in spacetime and their use for interstellar travel, A tool for teaching general relativity, Am. J. Phys. 56 (1988) 395-412.
\bibitem{LLO} P.S. Lemos, S.N. Lobo, Q. Oliveira, Morris-Thorne wormholes with a cosmological constant, Phys. Rev. D. 68 (2003) 064004.
\bibitem{T1} M. Thomas, Exact geometric optics in a Morris-Thorne wormhole spacetime, Phys. Rev. D. 77 (2008) 044043.
\bibitem{J} K. Jusufi, Conical Morris-Thorne wormholes with a global monopole charge, Phys. Rev. D. 98 (2018) 044016.
\bibitem{S1} Z.I. Szab\'{o}, Structure theorems on Riemannian spaces satisfying $R(X, Y)R=0,$ I. The local version, J. Differ. Geom. 17 (1982) 531-582.
\bibitem{S2} Z.I. Szab\'{o}, Classification and construction of complete hypersurfaces satisfying $R(X, Y)R=0,$ Acta Sci. Math. 47 (1984) 321-348.
\bibitem{S3} Z.I. Szab\'{o}, Structure theorems on Riemannian spaces satisfying $R(X, Y)R=0,$ II. The global version, Geom. Dedic. 19 (1985) 65-108.
\bibitem{AD} A. Adam\'{o}w, R. Deszcz, On totally umbilical submanifolds of some class of Riemannian manifolds, Demonstr. Math. 16 (1983) 39-59.
\bibitem{C} M.C. Chaki, On generalized quasi-Einstein manifolds, Publ. Math. Debr. 58 (2001) 683-691.
\bibitem{R1} H.S. Ruse, On simply harmonic spaces, J. Lond. Math. Soc. 21 (1946) 243-247.
\bibitem{R2} H.S. Ruse, On simply harmonic `kappa spaces' of four dimensions, Proc. Lond. Math. Soc. 50 (1949) 317-329.
\bibitem{R3} H.S. Ruse, Three dimensional spaces of recurrent curvature, Proc. Lond. Math. Soc. 50 (1949) 438-446.
\bibitem{SDM} E. Sabina, C. Dhyanesh, S. Mousumi, Curvature properties of Morris-Thorne wormhole metric, J. Geom. Phys. 174 (2022) 104457.
\bibitem{AC1} N.S. Agashe, M.R. Chafle, A semi-symmetric non-metric connection on a Riemannian manifold, Indian J. Pure Appl. Math. 23 (1992) 399-409.
\bibitem{AC2} N.S. Agashe, M.R. Chafle, On submanifolds of a Riemannian manifold with a semi-symmetric non-metric connection, Tensor 55 (1994) 120-130.
\bibitem{Wang1} Y. Wang, Affine connections of non-integrable distributions, Int. J. Geom. Methods Mod. Phys. 17 (2020) 2050127.
\bibitem{Wang2} Y. Wang, Curvature of multiply warped products with an affine connection, Bull. Korean Math. Soc. 50 (2013) 1567-1586.
\bibitem{G1} M. G{\l}ogowska, Semi-Riemannian manifolds whose Weyl tensor is a Kulkarni-Nomizu square, Publ. Inst. Math. 72 (2002) 95-106.
\bibitem{G2} M. G{\l}ogowska, On quasi-Einstein Cartan type hypersurfaces, J. Geom. Phys. 58 (2008) 599-614.
\bibitem{SRK} A.A. Shaikh, I. Roy, H. Kundu, On some generalized recurrent manifolds, Bull. Iranian Math. Soc. 43 (2017) 1209-1225.
\bibitem{DDHKS} F. Defever, R. Deszcz, M. Hotlo\'{s}, M. Kucharski, Z. Sent\"{u}rk, Generalisations of Robertson-Walker spaces, Ann. Univ. Sci. Budapest. E\"{o}tv\"{o}s Sect. Math. 43 (2000) 13-24.
\bibitem{SH1} A.A. Shaikh, H. Kundu, On warped product generalized Roter type manifolds, Balkan J. Geom. Appl. 21 (2016) 82-95.
\bibitem{SH2} A.A. Shaikh, H. Kundu, On generalized Roter type manifolds, Kragujevac J. Math. 43 (2019) 471-493.
\bibitem{K} D. Kowalczyk, On the Reissner-Nordstr\"{o}m-de Sitter type spacetimes, Tsukuba J. Math. 30 (2006) 363-381.
\bibitem{SC} A.A. Shaikh, D. Chakraborty, Curvature properties of Kantowski-Sachs metric, J. Geom. Phys. 160 (2021) 103970.
\bibitem{DG} R. Deszcz, M. G{\l}ogowska, Some examples of nonsemisymmetric Ricci-semisymmetric hypersurfaces, Colloq. Math. 94 (2002) 87-101.
\bibitem{SH3} A.A. Shaikh, H. Kundu, On equivalency of various geometric structures, J. Geom. 105 (2014) 139-165.
\bibitem{T2} S. Tachibana, A theorem on Riemannian manifolds of positive curvature operator, Proc. Jpn. Acad. 50 (1974) 301-302.
\bibitem{DGPSS} R. Deszcz, M. G{\l}ogowska, M. Plaue, K. Sawicz, M. Scherfner, On hypersurfaces in space forms satisfying particular curvature conditions of Tachibana
type, Kragujev. J. Math. 35 (2011) 223-247.
\bibitem{R4} R. Deszcz, On pseudosymmetric spaces, Bull. Belg. Math. Soc., Ser. A 44 (1992) 1-34.
\bibitem{R5} R. Deszcz, Curvature properties of a pseudosymmetric manifolds, Colloq. Math. 62 (1993) 139-147.
\bibitem{SH4} A.A. Shaikh, H. Kundu, On some curvature restricted geometric structures for projective curvature tensor, Int. J. Geom. Methods Mod. Phys. 15 (2018)
1850157.
\bibitem{SH5} A.A. Shaikh, H. Kundu, On warped product manifolds satisfying some pseudosymmetric type conditions, Differ. Geom. Dyn. Syst. 19 (2017) 119-135.
\bibitem{S} A.A. Shaikh, On pseudo quasi-Einstein manifolds, Period. Math. Hung. 59 (2009) 119-146.
\bibitem{SYH} A.A. Shaikh, D.W. Yoon, S.K. Hui, On quasi-Einstein spacetimes, Tsukuba J. Math. 33 (2009) 305-326.
\bibitem{B} A.L. Besse, Einstein manifolds, Springer-Verlag, Berlin, Heidelberg, 1987.
\bibitem{D} R. Deszcz, On some Akivis-Goldberg type metrics, Publ. Inst. Math. 74 (2003) 71-84.
\bibitem{DGJZ} R. Deszcz, M. G{\l}ogowska, J. Je{\l}owicki, Z. Zafindratafa, Curvature properties of some class of warped product manifolds, Int. J. Geom. Methods Mod.
Phys. 13 (2016) 1550135.
\bibitem{SH6} A.A. Shaikh, H. Kundu, On curvature properties of Som-Raychaudhuri spacetime, J. Geom. 108 (2016) 501-515.
\bibitem{G} A. Gray, Einstein-like manifolds which are not Einstein, Geom. Dedic. 7 (1978) 259-280.
\bibitem{SJ} A.A. Shaikh, S.K. Jana, On weakly cyclic Ricci symmetric manifolds, Ann. Pol. Math. 89 (2006) 139-146.
\bibitem{MM1} C.A. Mantica, L.G. Molinari, Extended Derdzinski-Shen theorem for curvature tensors, Colloq. Math. 128 (2012) 1-6.
\bibitem{MM2} C.A. Mantica, L.G. Molinari, Riemann compatible tensors, Colloq. Math. 128 (2012) 197-210.
\bibitem{MM3} C.A. Mantica, L.G. Molinari, Weyl compatible tensors, Int. J. Geom. Methods Mod. Phys. 11 (2014) 1450070.
\bibitem{MS1} C.A. Mantica, Y.J. Suh, The closedness of some generalized curvature $2$-forms on a Riemannian manifold I, Publ. Math. 81 (2012) 313-326.
\bibitem{MS2} C.A. Mantica, Y.J. Suh, The closedness of some generalized curvature $2$-forms on a Riemannian manifold II, Publ. Math. 82 (2013) 163-182.


\end{thebibliography}
\end{document}